\documentclass{amsart}
\usepackage{amssymb, array, verbatim, amscd}
\usepackage{amsmath, amsfonts,amsthm}

\newtheorem{theorem}{Theorem}[section]

\newtheorem{proposition}[theorem]{Proposition}

\newtheorem{definition}[theorem]{Definition}
\newtheorem{example}[theorem]{Example}

\newcommand{\Z}{{\mathbb Z}}
\begin{document}

\title[Application of Finite Fields]{Applications of Finite Fields to Dynamical Systems and
Reverse Engineering Problems }

\author{Mar\'{\i}a A. Avi\~n\'o }
\address{ Department of Mathematic-Physics\\
       University of Puerto Rico \\
       Cayey, PR 00736}
\email{mavino@cayey.upr.edu}
\author{ Edward Green}
\address{Mathematics Department\\
       Virginia Tech, Blacksburg,\\
       VA 24061-0123}
\email{green@math.vt.edu}
\author {Oscar Moreno}
\address{Department of Mathematics and Computer Science\\
University of Puerto Rico at Rio Piedras\\
Rio Piedras, PR 00931}
\email{moreno@uprr.pr}

\begin{abstract}
We present a mathematical model: dynamical systems over finite
sets (DSF), and we show that Boolean and discrete genetic models
are special cases of DFS.  In this paper, we prove that a function
defined over finite sets with different number of elements can be
represented as a polynomial function over a finite field. Given
the data of a function defined over different finite sets, we
describe an algorithm to obtain all the polynomial functions
associated to this data.  As a consequence, all the functions
defined in a regulatory network can be represented as a polynomial
function in one variable or in several variables over a finite
field. We apply these results to study  the reverse engineering
problem.
\end{abstract}

\subjclass{Primary:11T99 ; 05C20 }

\keywords{Finite fields, dynamical systems, partially defined
functions, regulatory networks}
\maketitle

\section{Introduction}
In this paper we introduce the definition of dynamical systems
over different finite sets (DSF) and we develop its applications
to \emph{regulatory networks} and the \emph{Reverse Engineering
Problem}.  We consider variables over sets with different numbers
of elements and we change that to variables over a finite field.

 The justification for considering dynamical systems over different finite
 sets is related with the method \emph{Generalized Logical Networks}
 developed by Thomas and colleagues, \cite{T,T1,T2,T3,T4}.
The  generalized logical networks has  a mean consideration: a
variable can have more that two possibilities but always the
number of possibilities is finite. In addition, the network is
described by a function which acts over several variables and for
each variable there are different number of values. These
considerations are very important for biologists because it is
known that in a regulatory network all the variables do not have
the same number of states. Here, we prove that all of these
functions can be considered over a finite field and as a
consequence of that we can represent them by polynomial functions.
In  section 2 we present an algorithm which changes a function
over different set of values to a function over a finite field.

In Section 2 we introduce the method to construct functions over a
finite field using functions defined over finite sets with
different number of elements. In Section 3, we apply the partially
defined functions to Reverse Engineering Problem. In Section 5, we
introduce the definition of Dynamical Systems over different
finite sets.
\section{Partially defined functions}
Now, in this section we introduce the mathematical background
which will permit the application of modelling  methods such as
generalized logical networks.

Let $X_j=\{0,1,\ldots ,j-1\}$ and let $\Z_p$ be the set of
integers modulo $p$, with $p$ a prime number. Suppose that $ p\geq
j$,   and we consider a canonical map from  $X_j$ to the field
$\Z_p$ given by $a \rightarrow a (\mod$ $ p)$. In the following we
consider $X_j\subset \Z_p$. Let $\textit{\textbf{x}}=(x_1,\ldots
,x_n)\in {\Z_p}^n$. We denote the polynomial ring in $n$ variables
over $\Z_p$ by $\Z_p[x_1,\ldots x_n]$. We begin with some
definitions. Let $D\subseteq \Z_p$.
\begin{definition}\label{ef}
Let $S\subsetneq {\Z_p}^n$. Let $f:S \rightarrow D$ be a function.
We will call $f$ a \emph{partially defined function} over $\Z_p$.
\end{definition}
\begin{example}\label{example}
  Now, let $g:X_2\times \Z_3\rightarrow X_2$, given by the
following table: \[\begin{tabular}{|c|ccc|}
  \hline
 $g(x_1,x_2)$ & 0 & 1 & 2 \\
  \hline
  0 & 1 & 0 &  1\\
  1 & 0 & 0 &  1\\
\hline
\end{tabular}.\]
Then the table of values of the partially defined function
$\hat{g}$ is the following:
 \[\begin{tabular}{|c|ccc|}
  \hline
 $\hat{g}(x_1,x_2)$ & 0 & 1 & 2 \\
  \hline
  0 & 1 & 0 &  1\\
  1 & 0 & 0 &  1\\
  2 & * & *  &*  \\
  \hline
\end{tabular}.\]
\end{example}
Let $S=X_{m_1}\times \cdots \times X_{m_n}\subsetneq {\Z_p}^n$,
and let $D\subseteq \Z_p$. Since a partially defined function
$\hat {f}:S\rightarrow D$ is not a function from ${\Z_p}^n$ to
$\Z_p$, we are interested in solving the following problem:

 [DF($\Z_p$):] Let $S\subsetneq {\Z_p}^n$ and let
$f:S\rightarrow D$ be a function. We want a polynomial function
$P:{\Z_p}^n\rightarrow \Z_p$ such that $P(
\textit{\textbf{x}})=f(\textit{\textbf{x}}), \ \hbox{ for all  }
\textbf{x}\in S\subsetneq {\Z_p}^n.$

 A function $P$ associated
to $f$ will be called \emph{a polynomial function for $f$}. Now,
we prove that the problem DF($\Z_p$) can have more than one
solution.
\begin{proposition}\label{function}
 For each function $f:S\rightarrow
D$  there is a polynomial $P(\textit{\textbf{x}})\in
\Z_p[x_1,\ldots ,x_n]$, such that $f(\textbf{x})=P(\textbf{x})$
for all $\textbf{x}\in S$. The polynomial $P$ can be chosen with
degree less than or equal to $n(p-1)$ but in general, it is not
unique.
\end{proposition}
\begin{proof}  If $k$ is a finite field  and
$f:k^n\rightarrow k$ is a function then there exists a polynomial
$P$ in the variables $x_1,\ldots ,\ x_n$, with coefficient in $k$,
such that $f(\textit{\textbf{x}})=P(x_1,\ldots , x_n)$ for all
$(x_1,\ldots , x_n)\in k^n $, \cite{LL}.  But, in our case we do
not have a function from $k^n\rightarrow k$, so  we will prove
that the polynomial function exists associated to the partially
defined function but it is not unique.

We will show the idea using the example \ref{example}. In the
table of $\hat{f}$ we can complete the table in some way. Then
there exists a unique polynomial for this table. But we can
complete the table in many ways, so the polynomial function exists
but, it is not unique.
\end{proof}
As a consequence of Proposition \ref{function} we have an
algorithm which solves the problem DF$(\Z_p)$. Let $f:S
\rightarrow \Z_p$ be a function. Let  $m=|S|$ be the cardinality
of $S$ ($m=\sum_{i=1}^nm_i$ when $f:X_{m_1}\times \cdots \times
X_{m_n} \rightarrow \Z_p$). Now,  we write a polynomial $P$ in $n$
variables $x_1,\ \ldots ,\ x_n$. $P$ has degree less than or equal
to $p-1$ in each variable, so has degree less than or equal to
$n(p-1)$. We denote $P$ in the following form: $P(x_1,\ldots
,x_n)=\sum_{\alpha \in {\Z_p}^n}b_\alpha x^\alpha ,$ where
$\alpha=(\alpha _1,\ldots ,\alpha _n), \ x^\alpha=x^{\alpha
_1}\ldots x_n^{\alpha _n}.$ Now, we evaluate $P$ for all
$\textbf{a}\in S$ and  we obtain a system of $m$ linear equations
in the $p^n$ unknowns $b_\alpha$ which always has solutions. The
system is the following:
\[(I) \ \sum_{\alpha \in {\Z_p}^n} (\textbf{a})^\alpha
b_\alpha=f(\textbf{a})\ \hbox{ for all }\textbf{a}\in S.\]
 Solving the system using elementary row operations, we finally obtain all the solutions.
 In \cite{G2}, it is proved that the rank of this system is $m$.
 Then, there are  $b_{\beta_1},\dots ,b_{\beta _m}$  coefficients of the polynomial $P$
 whose are determined in term of the  free coefficients denoted by $b_{\gamma_1},\dots ,b_{\gamma _{p^n-m}}$.
 Now, let ${\Z_p}^{(p-1)}[x_1,\ldots ,x_n]$ be the subspace of $\Z_p[x_1,\ldots ,x_n]$ of all  polynomials
with maximum degree $p-1$ in each variable and coefficients in
$\Z_p$. So, we have the following theorem.
 \begin{theorem}\label{function2}
 All the polynomial solutions with degree $\leq n(p-1)$ of the problem
 DF$(\Z_p)$ are given by a particular solution $f_0(\textit{\textbf{x}})$
  of $(I)$ plus the subspace \[U=\{g\in{\Z_p}^{(p-1)}[x_1,\ldots ,x_n]|g(\textbf{a})=0,\forall \textbf{a}\in
 S\}\] of dimension $p^n-m$.
 \end{theorem}
\begin{proof}
We know by linear algebra that  all the  solutions of  $(I)$ are
given by
 \[f_0(\textit{\textbf{x}})+b_{\gamma_1}g_1(\textit{\textbf{x}})
 +\cdots +b_{\gamma _{p^n-m}}g_{p^n-m}(\textit{\textbf{x}})\]
where $b_{\gamma_1},\ldots,b_{\gamma _{p^n-m}}\in \Z_p$ and $g_1,
\ldots, g_{p^n-m}\in U$. Let $h_1$ and $h_2$ be two polynomial
solutions of $(I)$. Then $h_1-h_2\in U$, so the theorem holds.
\end{proof}
\section{Reverse Engineering Problem  \\ over finite sets}
Now, we connect the problem DF$(\Z_p)$ with the Reverse
Engineering Problem over $\Z_p$. The problem for partially defined
functions is equivalent to the following.

 [\emph{P}($\Z_p$):] Given $\textbf{a}_1,\ldots,\textbf{a}_m\in
{\Z_p}^n$, $\textbf{b}=(b_1,\ldots ,b_m)\in {\Z_p}^m$, with
$m<p^n$. Find a polynomial $P\in \Z_p[x_1,\ldots,x_n]$ such that
$P(\textbf{a}_j)=b_j$ for $j=1,\ldots ,m$.

The problems \emph{P}($\Z_p$) and DF$(\Z_p)$ are equivalent. In
fact,  we only need to take
$S=\{\textbf{a}_1,\ldots,\textbf{a}_m\}$ and
$\textbf{b}=(f(\textbf{a}_1),\ldots ,f(\textbf{a}_m))$.

 The problem \emph{P}($\Z_p$) was solved  by E. Green in \cite{G}. He
called the problem \emph{P}($\Z_p$) for $(\textbf{a}_1,\ldots
,\textbf{a}_m;\textbf{b})$ and he proved that if \emph{P}($\Z_p$)
has solutions then the Reverse Engineering Problem over $\Z_p$ has
solutions.

Now, we define the Reverse Engineering Problem  over sets with
different number of elements. Let $\{k_j\}$ be a family of $n$
finite sets where $|k_j|=m_j$. We denote by $k=k_1\times \cdots
\times k_n$. Let $\textbf{r}_1,\ \ldots ,\ \textbf{r}_{m+1} \in
k$. We assume that the vectors $\textbf{r}_j=(r_{j1},\ldots
,r_{jn})$ are obtained by experiments (like microarray) and we
assume that $\textbf{r}_j$ determines $\textbf{r}_{j+1}$. Then the
\emph{Reverse Engineering Problem} over $k$ is to find  a function
$F=(f_1,\ldots ,f_n): k\rightarrow k$ such that
$F(\textbf{r}_j)=\textbf{r}_{j+1}$ for $j=1, \ldots, m$. But, we
rewrite:

[(REP)] The \emph{Reverse Engineering Problem} over $k$ is to find
polynomial functions $f_s: k\rightarrow k_s$ such that
\[  \  f_s(\textbf{r}_j)=r_{s,j+1} \ \hbox{for }j=1, \ldots, m\
\hbox { and }s=1,\ldots,n\]

 Now, we prove that if we can solve DF$(\Z_p)$, we can solve (REP)
and use the same algorithm. In fact, if $\Z_p$ is the field such
that $p\geq m_j$ for all $j$, we take the partially defined
functions $\hat{f}_s$ over $\Z_p$ considering $k_s\subseteq \Z_p$.
Let $\textbf{r}=(r_1,\ldots ,r_n)\in k$. So, we consider
$S\subsetneq {\Z_p}^n$ and $\hat{f}_s(\textbf{r}_j)=r_{s,j+1}$ for
$j=1, \ldots, m$ and $s=1,\ldots,n$.

 We have proved the following proposition.
\begin{proposition}
The reverse engineering problem over set with  different number of
elements  has polynomial solutions by Proposition \ref{function},
and Theorem \ref{function2}.
\end{proposition}
\begin{definition}
The matrix $A=(\textbf{r}_j)_{m\times n}$ will be called  the
matrix of the problem  REP.
\end{definition}

\begin{example}
  Suppose we
have the following data: $\textbf{r}_1=(1,2,0)$,
$\textbf{r}_2=(2,2,1)$, $\textbf{r}_3=(1,0,1)$,
$\textbf{r}_4=(0,1,1)$, and  $\textbf{r}_5=(1,1,0)$. And, we have
the additional information:\\
(a) the variables $\{x,y,z\}$ are defined over different finite
sets, but we take  finite fields:  $x,y\in \Z_3= \{0,1,2\}$ and
$z\in \Z_2=\{0,1\}$.\\
(b)the variable $x$ depends of $x$ and $z$, $y$ depends of $x$ and
$y$, and $z$ depends of $y$ and $z$.
\end{example}
The matrix $A$ of the problem is the following:
\[ A=\left(\begin{array}{ccc}
1&2&0\\
2&2&1\\
1&0&1\\
0&1&1\\
1&1&0\cr
\end{array}\right)\]
We want polynomial functions $f_1,$ $f_2,$ and $f_3$, such that
$F=(f_1,f_2,f_3)$ and $F(\textbf{r}_j)=\textbf{r}_{j+1}$ for
$j=1,2,3,4$.

 The additional information (b)  means that the functions that we are looking
 for are as follows:
\[\begin{array}{ccc}
f_1(x,z)&=& a_0+a_1x +a_2z+
a_3xz+a_4x^2+a_5z^2\\
&& +a_6x^2z+a_7xz^2+a_8x^2z^2 \\
f_2(x,y)&= & b_0+b_1x +b_2y+
b_3xy+b_4x^2+b_5y^2\\
& & +b_6x^2y+b_7xy^2+b_8x^2y^2\\
f_3(y,z)& =&  c_0+c_1z
+c_2y+c_3yz+c_4y^2+c_5z^2\\
& & +c_6y^2z+c_7yz^2+c_8y^2z^2 \end{array} \]
 Using the data we
have the table of $f_1$.
\[\begin{tabular}{|c|ccc|}
\hline
$f_1(x,z)$ & 0 & 1 &2 \\
  \hline
  0 & * & 1 & * \\
 1 & 2 & 0 & * \\
  2 & *& 1 & * \\
  \hline
\end{tabular} \]
Using the above table and the algorithm for problem DF($\Z_p$), we
 obtain a system of $4$ linear equation with $9$ unknown. The  matrix of the system of linear
 equation is:
\[ \overline{A_1}=\left(\begin{array}{ccccccccc|c}
1&1&0&0&1&0&0&0&0&2\\
1&2&1&2&1&1&1&2&1&1\\
1&1&1&1&1&1&1&1&1&0\\
1&0&1&0&0&1&0&0&0&1\cr
\end{array}\right)\]
 Using elementary row operations, we have the following:
 \[ a_2 = 1+2a_3+2a_5+2a_6+2a_7+2a_8,\ a_4 = 1+2a_6+2a_8,\]
 \[  a_1 = 1+2a_3+2a_7, \ a_0 =
 a_3+a_7+a_6+a_8.\]
 A particular  solution of the system is $f_1=x+z+x^2$. And all the solutions
 are given by $f_1+a_3 g_1+a_5g_2+a_6g_3+a_7g_4+a_8g_5 $,
 where
\[ g_1=1+2x+2z+xz, \  g_2=2z+z^2, \  g_3=1+2z+2x^2+x^2z \]
 \[ g_4=1+2x+2z+xz^2, \hbox{ and }g_5=1+2z+2x^2+x^2z^2. \]
If we denote by $U_1$ the subspace of $\Z_3[x,y,z]$ generated by
$\{g_1,g_2,$ $ g_3,g_4,g_5\}$, then all the solutions with degree
$\leq 2$ in each variable, are $f_1+U_1$.

 Similarly we obtain:

(1)  $f_2=x+y^2$ and all the solutions are
 $f_2+U_2$, where $U_2$ is the subspace generated by the
 polynomials
 \[ h_1=2+x+xy+y^2,\ h_2=2+2y+x^2+2y^2,\]
  \[ h_3=1+2x+2y^2+xy^2, \
 h_4=y+2y^2+x^2y, \ h_5=2y+y^2+x^2y^2 . \]

(2)  $f_3=1+y+y^2$ and all the solutions are
 $f_3+U_3$, where $U_3$ is the subspace generated by the
 polynomials
 \[ v_1=2+z+2y+yz, \ v_2=1+2z+2y^2+y^2z,\]
 \[v_3=2+z+2y+yz^2,\ v_4=1+2z+2y^2+y^2z^2, \ v_5=2z+z^2 . \]

Finally one of the functions that can describe the genetic network
is the following\[f(x,y,z)=(x+z+x^2,x+y^2,1+y+y^2)\]

\section{Solution over the finite field GF($p^n$)}
We can solve the problems DF$(\Z_p)$ and (REP) using Lagrange
interpolation over the field GF$(p^n)=K$ \cite{M,O}. Let
$f:S\rightarrow \Z_p$ be a function with $S\subsetneq {\Z_p}^n$.
Let $|S|=m$ be the cardinality of $S$. let $\{ \alpha _1, \alpha
_2,\ldots ,\alpha _n\}$ be a fixed basis of $K$. There is a
natural one to one correspondence between the sets ${\Z_p}^n$ and
$K$, namely
\[\lambda:(a_1,\ldots,a_n)\mapsto a_1\alpha _1 +\cdots
+a_n\alpha _n.\]

Let $\overline {S}=\lambda (S)\subsetneq K$. Now we have the
partially defined function $\hat{f}=\lambda \circ f\circ \lambda
^{-1}:\overline {S} \rightarrow\Z_p $. We denote the elements of
$\overline {S}$ by $\overline {\textit{\textbf{a}}}$.

Now, using the Lagrange interpolation formula we have the
following: $\overline {\textit{\textbf{a}}}_1$, $\ldots $,
$\overline {\textit{\textbf{a}}} _m$ are $m$ distinct elements of
the finite field $K$ and $\hat{f}(\overline
{\textit{\textbf{a}}}_1)=b_1$, $\dots $, $\hat {f}(\overline
{\textit{\textbf{a}}}_m)=b_m$, with $b_1,\ldots ,b_m$ elements in
 $\Z_p$. We know that $\Z_p\subset K$. We rewrite the problem
DF$(\Z_p)$ as follows:

[DF($p^n$):] Let $\overline {S}\subsetneq K$ and let $f:\overline
{S}\rightarrow K$ be a function. We want a polynomial $P(x)\in
K[x]$ such that
\[P( \overline{\textit{\textbf{x}}})=f(\overline{\textit{\textbf{x}}}), \ \hbox{ for all  }
\overline {\textit{\textbf{x}}}\in \overline {S}\subsetneq K.\]

We can observe that this is the same problem (REP) if we consider
$\overline{S}=\{\textbf{r}_1,\ldots, \textbf{r}_{m+1}\}$ and
$f(\textbf{r}_j)=\textbf{r}_{j+1}\in K$. So in the following we
denote both problem by DF$(p^n)$.

 Using Lagrange Interpolation, we know that: there exists a
 polynomial $\overline {P}\in K[x]$ of degree $d\leq m-1$ such
 that $\overline {P}(\overline {\textit{\textbf{a}}}_i)=\textit{\textbf{b}}_i\in K$ for $i=1,
 \ldots ,m$. The polynomial is given by \[\overline{P}(x)=\sum
 _{i=1}^{m-1}\textit{\textbf{b}}_i\prod _{k=1,k\ne i}^{m-1}(\overline
 {\textit{\textbf{a}}}_i-\overline {\textit{\textbf{a}}}_k)^{-1}(x-\overline
 {\textit{\textbf{a}}}_k).\]
 Then, we have proved the following theorem.
 \begin{theorem}\label{LI}
 The problem DF$(p^n)$ has solutions over the field $K=$GF$(p^n)$
 using Lagrange Interpolation. That is, there exists a polynomial
 $\overline{P}_0\in K[x]$, such that
 $\overline{P}_0(\overline{\textbf{a}}_i)=\textit{\textbf{b}}_i\in K$, for $i=1,\ldots
 ,m$. The degree of the polynomial $\overline{P}_0$ is less than or
 equal to $m-1$. If  $I$ is the ideal of $k[x]$ generated by $\overline
 {P}(x)=(x-\overline
 {\textit{\textbf{a}}}_1)\cdots (x-\overline
 {\textbf{a}}_m)$, then all the solutions are given by
 $\overline{P}_0(x)+G(x)$, where $G(x)\in
 I$.
 \end{theorem}
 Now, we have a new algorithm to solve the problem. We know by the
 Theorem \ref{LI} that the solution is a polynomial of the
 following form: $\overline{P}(x)=\sum_{k=0}^{m-1}B_kx^k$.
We evaluate $\overline{P}(x)$
 in all the elements of $\overline{S}$. Then we obtain a system
of linear equations with one solution  \[(II)\ \hbox{   } \sum
_{k=0}^{m-1}B_k\overline
 {\textit{\textbf{a}}}_i^k=b_i\ \hbox {  for  } i=1,\ldots ,m\]
We want to remark that the system $(II)$ has rank $m$, since
$\overline {\textit{\textbf{a}}}_i\ne \overline
 {\textit{\textbf{a}}}_j$ for $i\ne j$.
 Finally we have an output a polynomial in one variable with
 degree less than or equal to $m-1$.
 \begin{example} Let
 $S=\{(0,1),(1,0),(1,1),(2,1)\}\subsetneq {\Z_3}^2$, and $ f:S\rightarrow\Z_3$. The function $f$ has  the following table of
 values.
 \[\begin{tabular}{|c|ccc|}
\hline
$f(x_1,x_2)$ & 0 & 1 &2 \\
  \hline
  0 & * & 1 & * \\
 1 & 2 & 0 & * \\
  2 & *& 1 & * \\
  \hline
\end{tabular} \]
\end{example}
Let $\alpha$ be  a root of the polynomial $X^2+X+2$ in $\Z_3$.
Then $\alpha ^2=2\alpha +1$, and a basis for $GF(3^2)$ is
$\{\alpha , 1 \}$. A natural correspondence is $(x_1,x_2)\mapsto
x_1 \alpha +x_2$ We have $\textbf{a}_1=(0,1)$,
$\overline{\textit{\textbf{a}}}_1=1$;
$\textit{\textbf{a}}_2=(1,0)$,
$\overline{\textit{\textbf{a}}}_2=\alpha$;
$\textit{\textbf{a}}_3=(1,1)$,
$\overline{\textit{\textbf{a}}}_3=\alpha +1$; and
$\textit{\textbf{a}}_4=(2,1)$,
$\overline{\textit{\textbf{a}}}_4=2\alpha +1$. Then, the
particular solution is given by a polynomial
$P_0(x)=c_0+c_1x+c_2x^2+c_3x^3$. We evaluate in the four elements
of GF$(3^2)$ using the table of values. The matrix of the system
is the following: \[ M=\left(\begin{array}{cccc|c}
1&1&1&1&1\\
1&\alpha&\alpha ^2 &\alpha ^3 & 2 \\
1&\alpha +1&2\alpha +1& 2\alpha &0\\
 1&2\alpha +1&2& \alpha +2& 1\cr
\end{array}\right)\]
Solving the system we obtain the polynomial $ P_0(x)=\alpha ^3
+2x+ \alpha ^6x^2 +x^3$. Now, we change the polynomial in two
polynomials with two variables $x_1,$ $x_2$ and coefficients in
$\Z_3$. We use the correspondence $\lambda^{-1}$ and obtain the
following:
\[f(x_1,x_2)=(2+x_1+2x_1x_2+x_2^2,2+2x_1+x_1^2+2x_1x_2+2x_2^2).\]

\section{Dynamical Systems over  finite sets with different number of elements}
In this section we present a  definition of dynamical system over
finite sets.
\begin{definition}
A \emph{ dynamical system over finite sets} is a pair $(S,f)$ such
that\\
(a) $S\subseteq {\Z_p}^n$.\\
(b) A function, $f=(f_1,\ldots , f_n):S\rightarrow S,$
 and each function $f_i:S_i\rightarrow \Z_p$ is a
partially defined function with $S_i\subseteq S$ for all $i$.
\end{definition}
  The DSF is a time discrete dynamical system, that is the dynamics is generated by
iteration of the function $f$.
\begin{definition}
The \textit{state space} $\mathcal{S}_{f}$ of the dynamical system
$f:S\rightarrow S$ is a finite directed graph (digraph) with
vertex set $S$ and arc set $A=\{(x,y)\in S:f(x)=y\}$, that is the
ordered pair $(x,y)$ stands for and arrow from a vertex
$x=(x_1,\ldots ,x_n)\in S$ to a vertex $y=(y_1,\ldots ,y_n)\in S$
if and only if $f(x)=y$.
\end{definition}
\begin{theorem}
Let $(S,f)$ be a dynamical system over finite sets. Then

(a) $f$ can be represented as a polynomial function in one
variable over the finite field GF$(p^n)$, for some  prime $p$.

(b) If $f=(f_1,\ldots,f_n)$ then each function $f_i$ can be
represented as a polynomial function in one variable or in several
variables over $\Z_p$, for some  prime $p$.

(c) Part $(a)$ and $(b)$ hold for all prime number $p\geq max_i
|\overline{S_i}|$, where $|\overline{S_i}|$ is the maximal number
of different elements in the  coordinate $i$ of the set $S_i$
\end{theorem}
\begin{proof}It is a consequence of the Sections 1  and 2.
\end{proof}

\section{Examples and  applications}
In this section we present two applications of the representation
of the dynamical systems by polynomials over a finite fields. One
very important things is to determine the steady states, that is
the elements $\textbf{x}$ such that $f(\textbf{x})=\textbf{x}$. We
determine that in the first example.

In the reverse engineering problem, we have a set of solutions and
here we suggest a method for biologist to determine if one of the
solution is the right one.

 In Fig. 1, an example of regulatory network is shown. This example of Generalized
   Logical Networks appear in \cite {J}. Here,  we  use the usual words for biologists. Gene $1$
  regulates genes $2$ and $3$, so that it has
two \emph{thresholds} (two values different $0$) and the
corresponding \emph{logical variable} $x_1$ takes its value from
$\{0, 1, 2\}$. Similarly, $x_2$ and $x_3$ have one and two
thresholds, respectively, and hence possible values $\{0,1\}$ and
$\{0, 1, 2\}$.   The functions $f_1({\bf x})$, $f_2({\bf x})$, and
$f_3({\bf x})$ need to be specified such as to be consistent with
the \emph{threshold restrictions} in the graph. Examples of
\emph{logical functions} allowed by the generalized logical method
are shown in Fig. 1(b).

 Consider the case of $f_2({\bf x})$. If
$x_1 \ne 0$ and $x_3 \ne  0$, so that $x_1$ and $x_3$ have values
\emph{above their first threshold}, the \emph{inhibitory
influences of genes } $1$ and $3$ on gene $2$ become
\emph{operative}. Figure 1(b) indicates that $x_2$ will
\emph{tend} to $0$, that is, below the first threshold of the
protein produced by gene $2$. If either $x_1= 0$ or $x_3= 0$, that
is, if only one of the inhibitory influences is operative, then
gene $2$ is \emph{moderately expressed}. This is here represented
by the value $1$ for the image of $x_2$. In general, several
logical functions will be consistent with the threshold
restrictions. Exactly which logical function is chosen may be
motivated by biological considerations or may be a guess
reflecting uncertainty about the structure of the system being
studied.
\begin{example}
 \[\ \ \ \ \
\begin{array}{c}
G_3 \circlearrowleft  \\
\nearrow     \searrow  \\
G_1 \rightleftarrows  G_2\\
 \hbox {Figure 1 (a)} \cr
\end{array}  \hbox{     } \hbox{     }
\begin{tabular}{|c|c|c|}
  \hline
  $f_1(x_2)$ & 0 & 2 \\
  \hline
  $x_2$& 0 & 1 \\
  \hline
\end{tabular} \]
\[ \begin{tabular}{|c|ccc|}
\hline
$f_2(x_1,x_3)$ & 0 & 1 & 2 \\
  \hline
  0 & 1 & 1 & 1 \\
1 & 1 & 0 & 0 \\
 2 & 1 & 0 & 0 \\
  \hline
\end{tabular} \\ \hbox{     } \hbox{     }
 \begin{tabular}{|c|ccc|} \hline

 $f_3(x_1,x_3)$ & 0 & 1 & 2 \\
  \hline
  0 & 2 & 2 & 1 \\
  1 & 2 & 2 & 1 \\
  2 & 0 & 0 & 0 \\
  \hline
\end{tabular}\]
\centerline{Figure 1 (b)}
\end{example}
Now, we describe this  example   using a dynamical system over
finite sets. We have three genes, and the regulatory network is
the following:

(1) The digraph: $Y \ \  \ \ \ \begin{array}{ccc}
& & 3 \circlearrowleft\\
 & \nearrow & \downarrow \\
1 & \rightleftarrows & 2
\end{array}$  \\
(2) The variables $x_1$ and $x_3$ are in  $\Z_3=\{0,1,2\}$, and
the variable $x_2$ is in  $X_2=\{0,1\}$,\\
(3) There are several possibilities for functions $f_1$ and $f_2$,
but $f_3$ is the unique function. Then  the polynomial functions
$f_1$, $f_2,$ and $f_3$ are:
\[\begin{array}{ccl}
 f_1(x_2)&= &2x_2  \\
 f_2(x_1,x_3)&= &1+2x_1^2x_3^2 \\
f_3(x_1,x_3)&=
&2+x_1+2x_3+x_1x_3+2x_1^2+x_3^2+2x_1^2x_3\\
& & +2x_1x_3^2+x_1^2x_3^2
\end{array}\]
(4) The global function $f:(\Z_3)^3 \rightarrow (\Z_3)^3$,
\[f(x_1,x_2,x_3)= (f_1(x_2),f_2(x_1,x_3),f_3(x_1,x_3)).\]
The state space has vertices $V=\Z_3\times X_2\times \Z_3$. We
want to know the steady states of the dynamical system $f$, in
general it is a very difficult problem. But, in this particular
case  we have three equations in three variables,
\[2x_2=x_1, \ 1+2x_1^2x_3^2=x_2, \]
\[2+x_1+2x_3+x_1x_3+2x_1^2+x_3^2+2x_1^2x_3+2x_1x_3^2+x_1^2x_3^2=x_3\]
For $x_2$ we have only two values,  if $x_2=0$ then $x_1=0$ and we
obtain $0=1$ in the second equation, that is impossible. If
$x_2=1$ then $x_1=2$ and $1+2x_3^2=1$ so $x_3=0$. We can check in
the last equation and the only solution is $(2,1,0)$.
\begin{example} In the example that we present in Section 3, we
have $3^{15}=14,348,907$ different solutions. But, we can select
the particular solution $f(x,y,z)=(x+z+x^2,x+y^2,1+y+y^2)$, and
try to find which vectors in the state space go to the first state
$(1,2,0)$. \end{example} Solving over $\Z_3$ the equations
$x+z+x^2=1,$ $x+y^2=2,$ and  $1+y+y^2=0$, we obtain that:
$1+y+y^2=0$ has one solution $y=1$. So $x=1$ and $z=2\equiv 0
(\mod 2)$. Therefore, $f(1,1,0)=(1,2,0)$. Now, we can check that
in the laboratory. Since the state space has only $18$ elements,
there are several functions with this property.

%\section{Acknowledgments}


\begin{thebibliography}{99}
\bibitem{A}T. Akutsu, S. Kuahara, O. Maruyama, and S. Miyano,
\textit{Identification of gene regulatory networks by strategic
gene disruptions and gene overexpressions.} Proceedings of the 9th
ACM-SIAM Symposium on Discrete Algorithms, H. Karloff, ed. (ACM
Press, 1998).

\bibitem{G} Edward. Green,\textit{On polynomial solutions to
reverse engineering problems.} to appear.

\bibitem{G2} E. L. Green, M. Leamer, and A. Li,\textit{Polynomial models of times series over $(\Z/p)^n$,}
to appear.

\bibitem{K} T.E. Ideker, V. Thorsson, and R.M. Karp.\textit{ Discovery
of Regulatory Interactions Through Perturbation: Inference and
Experimental Design,} Pacific Symposium on Biocomputing 5:302-313
(2000).

\bibitem{J} H. de Jong,\textit {Modeling and Simulation of Genetic
Regulatory Systems: A Literature Review} INRIA, No. 4032, 2002.

\bibitem{LS2} R. Laubenbacher, J. Shah, and B. Stigler,\textit{ Simulation of
polynomial systems.} To appear in {\em Simulation in the Health
and Medical Sciences,} Society for Modeling and Simulation
International, San Diego, CA, 2003.

\bibitem{LL} R. Lidl and H. Niederreiter, \textit{ Introduction to
finite fields and their applications}. Cambridge University Press.
1994.

\bibitem{M} O. Moreno, D. Bollman, and  M. Avi\~no-Diaz, \textit{ Finite
Dynamical Systems, Linear Automata, and Finite Fields,} Conference
Proceedings, WSEAS Transactions, (2002).

\bibitem{O} H. Otiz-Zuzuaga, M. Avi\~no-Diaz, C. Corrada,
R. Laubenbacher, S. Pe\~na de Ortiz,  and O. Moreno,
\textit{Applications of finite fields to the study of microarray
expression data} preprint(2003).

\bibitem{S}R. Somogy and C. Sniegoski,\textit{ Modeling the complexity
of genetic networks: understanding multigenic and pleiotropic
regulation.} Complexity (1): 45-63, 1996.

\bibitem{T} R. Thomas, \textit{Regulatory networks seen as asynchronous automata: A
logical description}. J. Theor. Biol. 153, 1–23,  1991.

\bibitem{T5} R. Thomas, \textit{Laws for the dynamics of regulatory
networks.} Int. J. Dev. Biol. 42, 479–485, 1998.

\bibitem{T1} R. Thomas, and R. d'Ari,
\textit{Biological Feedback,} CRC Press, Boca Raton, FL,1990.

\bibitem{T2} R. Thomas, A.-M. Gathoye, and L. Lambert\textit{ A complex
control circuit: Regulation of immunity in temperate
bacteriophages.} Eur. J. Biochem. 71, 211–227, 1976.

\bibitem{T3} R. Thomas, and M. Kaufman, \textit{ Multistationarity, the basis
of cell differentiation and memory: I. Structural conditions of
multistationarity and other nontrivial behavior. Chaos 11(1)},
170–179, 2001.

\bibitem{T4} R. Thomas, D. Thieffry, and M. Kaufman, \textit{ Dynamical
behaviour of biological regulatory networks: I. Biological rool of
feedback loops and practical use of the concept of the
loop-characteristic state.} Bull. Math. Biol. 57(2), 247– 276,
1995.


\end{thebibliography}
\end{document}